\documentclass[ prb, reprint]{revtex4-1}

\usepackage{graphicx}
\usepackage{subfigure}
\usepackage{amsmath,amssymb,amsthm}
\usepackage[usenames,dvipsnames,svgnames,table]{xcolor}
\usepackage{bm}
\usepackage{framed}
\usepackage[normalem]{ulem}
\usepackage{cancel}
\usepackage{mathtools}

\newcommand{\va}{\boldsymbol{\mathsf{a}}}

\newcommand{\ve}{\boldsymbol{e}}

\newcommand{\vq}{\boldsymbol{q}}

\newcommand{\vw}{\boldsymbol{w}}
\newcommand{\vx}{\boldsymbol{x}}
\newcommand{\vz}{\boldsymbol{z}}

\newcommand{\rr}{\mathsf{r}}
\newcommand{\qq}{\mathsf{q}}

\newcommand{\vB}{\boldsymbol{B}}

\newcommand{\vC}{\boldsymbol{C}}

\newcommand{\vF}{\boldsymbol{\mathsf{F}}}
\newcommand{\vFQ}{\boldsymbol{\mathsf{FQ}}}

\newcommand{\vR}{\boldsymbol{\mathsf{R}}}
\newcommand{\vT}{\boldsymbol{\mathsf{T}}}
\newcommand{\vS}{\boldsymbol{\mathsf{S}}}
\newcommand{\vQ}{\boldsymbol{\mathsf{Q}}}

\newcommand{\vX}{\boldsymbol{X}}

\newcommand{\QQ}{\mathbb{Q}}
\newcommand{\RR}{\mathbb{R}}
\newcommand{\CC}{\mathbb{C}}

\newcommand{\diag}{\text{diag}}

\theoremstyle{definition}

\theoremstyle{remark}

\begin{document}

\title{Rank-1 accelerated illumination recovery in scanning diffractive imaging by transparency estimation.}

\begin{abstract}
We consider the problem of blind ptychography, that is the joint estimation of an unknown object and an illumination function
from diffraction intensity measurements.  In ptychography, diffraction measurements from neighboring regions of the same object are related to each other 
by a pairwise relationship between overlapping frames. When the illumination is well known, the relationship among frames is given by a linear projection operator.
We propose a power iteration-projection algorithm that minimizes the global pairwise discrepancy among frames. We accelerate the convergence of power method by  subtracting the estimated localized average transparency of the unknown object. 
The method is effective for weakly scattering and low contrast objects or piecewise smooth specimens.
\end{abstract}

\author{Stefano Marchesini }
\affiliation{Advanced Light Source, Lawrence Berkeley National Laboratory, Berkeley, CA 94720 (USA)}
\author{Hau-tieng Wu}
\affiliation{Department of Mathematics, University of Toronto, Toronto, Ontario (CANADA)}
\maketitle

\begin{widetext}
Ptychography is an increasingly popular technique to achieve diffraction limited imaging over a large field of view without the need for high quality optics
\cite{Hoppe:1969,Hegerl_Hoppe:1970,batesptycho, spence:book,Chapman1996,Thibault2008,Thibault2009, pie2,Kewish2010,Honig:11,guizar-mirror,Maiden:2012, axel,Thibault-Guizar,Godard:12,clark:2011,fourierptycho,batey2014information, dong2014spectral, marrison2013ptychography, tian2014multiplexed,vine2009ptychographic, stockmar2013near}.
Since the reconstruction of ptychographic data is a non-linear problem, there are still many open problems\cite{synchroAP}, nevertheless the phase retrieval problem is made tractable by recording multiple diffraction patterns from the same region of the object, compensating phase-less information with a redundant set of measurements. Data redundancy enables to handle experimental uncertainties as well. 
Methods to work with unknown illuminations or ``lens'' were proposed\cite{Chapman1996,Thibault2008,Thibault2009, pie2, hesse2015proximal}. 
They are now used to calibrate high quality x-ray optics  \cite{Kewish2010,Honig:11,guizar-mirror}, EUV lithography tools\cite{wojdyla:2014}, x-ray lasers \cite{ptychoxfel}  and space telescopes \cite{spacetelescopes}.
More recently,  position errors\cite{fienup,Maiden:2012, axel}, background\cite{thurman2009,guizar:bias}, noise statistics \cite{Thibault-Guizar,Godard:12} and partially coherent illumination\cite{clark:2011,Fienup:93,Abbey2008,CDIpartialcoherence}. Situations when sample, illumination function, incoherent multiplexing effects, as well as positions,  vibrations, binning, multiplexing, fluctuating background are  unknown parameters in high dimensions  have been added to the nonlinear optimization to fit the data  using projections, gradient, conjugate gradient, Newton \cite{chao},   and spectral methods \cite{marchesini2013augmented,synchroAP,batey2014information} . 
Here we focus on the illumination retrieval problem. We utilize the notation described in \cite{marchesini2013augmented,synchroAP}, which is summarized below.

\begin{figure}[hbtp]
  \begin{center}
  	\includegraphics[width=.4\textwidth]{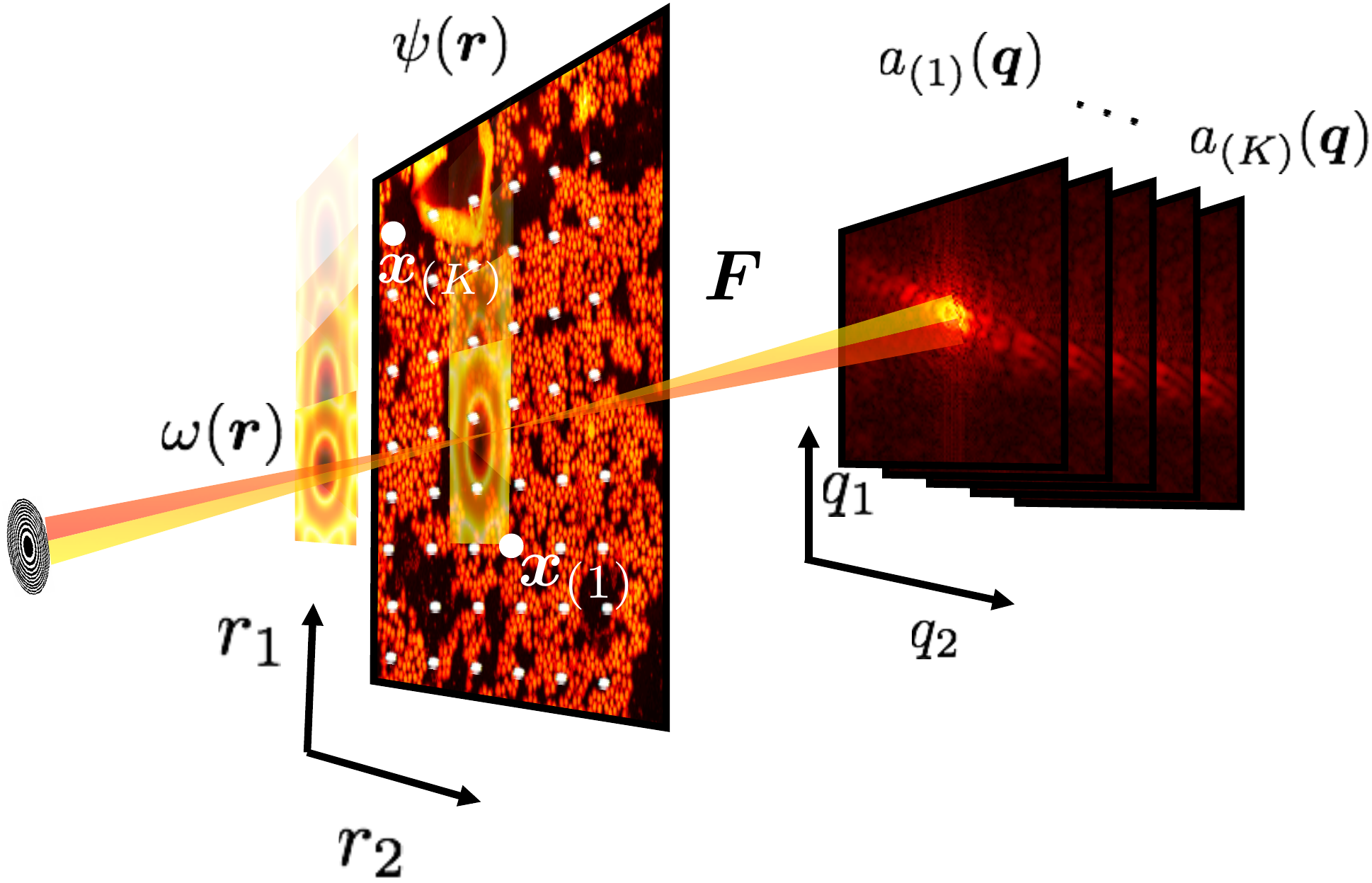}
          \includegraphics[width=0.4\textwidth]{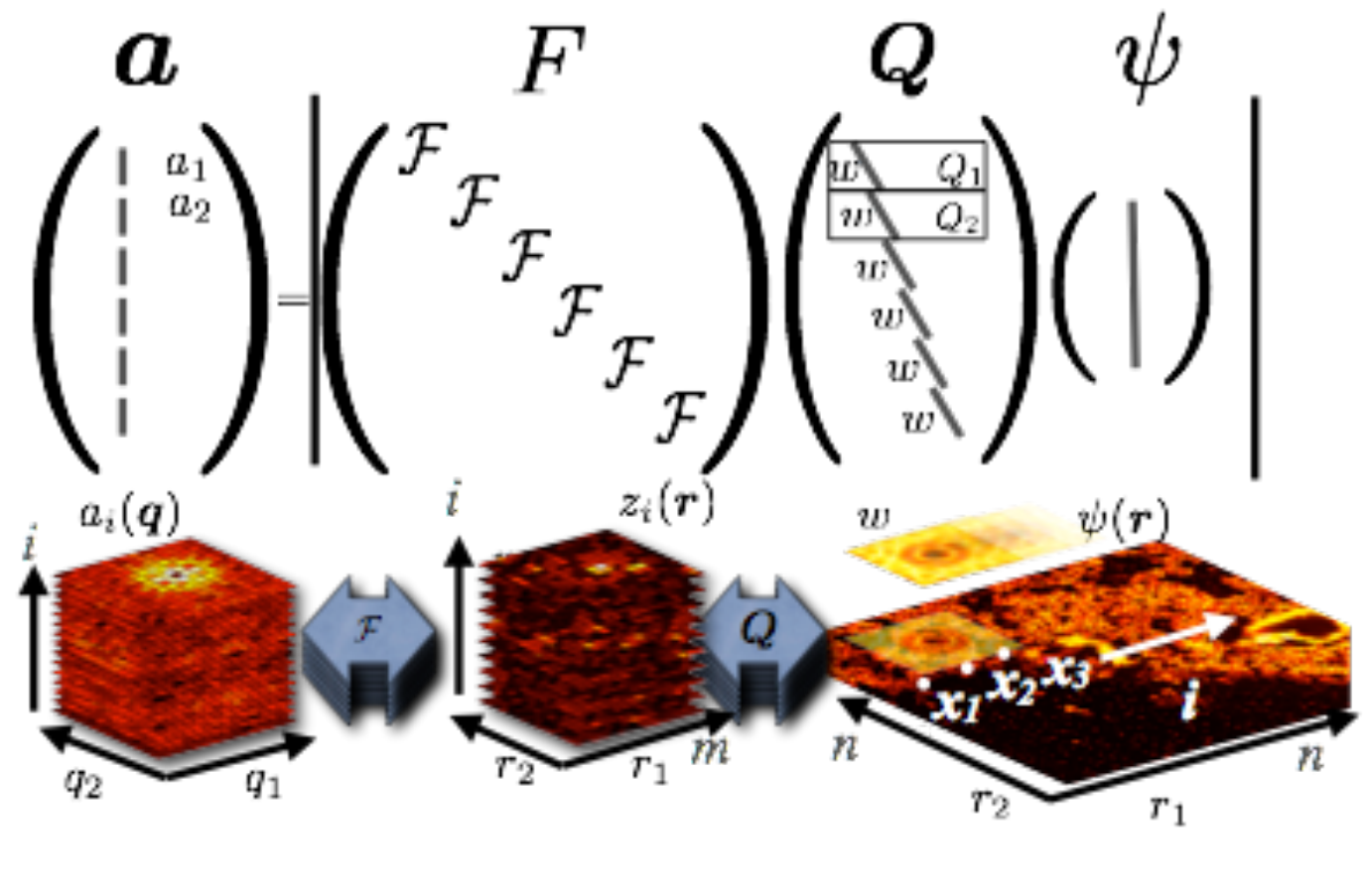}
 \end{center}
 \caption{
\label {fig:exp_eqs}
 The measured amplitudes $\va$ and the relationship with an an unknown object $\psi$ in standard ``Far Field" ptychography.
 Geometric representation of the  operators involved in simulating a ptychographic imaging experiment
 \cite{synchroAP,marchesini2013augmented}.
 }
\end{figure}

\section{Background}

The relationship between an unknown discretized object $\psi$ , and the diffraction measurements $\va_{(i)}$ collected in a ptychography experiment 
(see figure \ref{fig:exp_eqs}) can be represented compactly as:
\begin{align}
\label{eq:forwardproblem}
&\va = |\vFQ\psi^\vee| \text{, or } \begin{cases}
\va = |\vF\vz|,\\
\vz = \vQ\psi^\vee,
\end{cases}
\end{align}
where $\psi^\vee$ is the object $\psi$ in vector form.  Eq. (\ref{eq:forwardproblem}) can be expressed as:
\begin{align*}
&\overbracket{\left[\begin{array}{c} \va_{(1)}\\ \vdots \\ \va_{(K)}\end{array} \right]}^{\va\in \RR^{Km^2}}=\left |
\overbracket{
\left[\begin{array}{ccc} F & \ldots & 0 \\ \vdots & \ddots & \vdots\\ 0 & \ldots & F\end{array} \right]
}^{\vF\in \CC^{Km^2\times Km^2}}
\overbracket{
\left[\begin{array}{c} \vz_{(1)}\\ \vdots \\ \vz_{(K)}\end{array} \right] }^{\vz \in \CC^{K m^2}}
\right |,
\quad
\overbracket{\left[\begin{array}{c} \vz_{(1)}\\ \vdots \\ \vz_{(K)}\end{array} \right] }^{\vz \in \CC^{K m^2}}
=\overbracket{ \left[\begin{array}{c} \diag(w) \vT_{(1)}\\ \vdots \\ \diag(w)  \vT_{(K)}\end{array} \right]}^{\vQ \in \CC^{K m^2\times n^2},}
\overbracket{\left[\begin{array}{c} \psi_{1}\\ \vdots \\ \psi_{n^2}\end{array} \right] }^{\psi \in \CC^{n^2}}
\end{align*}
where $\vz$ are $K$ frames extracted from the object $\psi$ and multiplied by the illumination function $w$, $\vF$ is the associated 2D DFT matrix when we write everything in the stacked form, that is, $F$ is a $m^2\times m^2$ matrix satisfying $F_{lk}=\ve_l^TF\ve_k=e^{i\qq_m^{-1}(l-1)\cdot \rr_m^{-1}(k-1)}$. 
For experimental geometries related to what has been described above, such as Near Field, ``Fresnel", Fourier ptychography, through-focus \cite{dong2014spectral, vine2009ptychographic, stockmar2013near,tian2014multiplexed} one can substitute the simple Fourier transform with the appropriate propagator.

The illumination matrix  $\vQ$ encodes information about  illumination $w\in \CC^{m^2}$, and the $K$ relative translations between the probe $w$ and the object $\psi$. In particular, the $(i)$  block denoted as $\vQ_{(i)}\in \CC^{m^2\times n^2}$ extracts one frame from the object and multiplies by the illumination $w$:
\begin{align}\label{definition:Qi:rewrite}
\vQ_{(i)}=&
\overbracket{{\diag}(w^\vee) }^{\boxed{\bcancel{\quad}} }
\overbracket{\vT_{(i)} \in \CC^{m^2\times n^2}}^{\boxed{\quad \bcancel{\quad}\quad \quad}}
,\quad
\vQ= 
\left (
\begin{array}{c}
 {\diag}(w^\vee)  \vT_{(1)}\\
{\diag}(w^\vee)   \vT_{(2)}\\
  \vdots\\
{\diag}(w^\vee) \vT_{(K)}
\end{array}
\right ) 
\quad
\end{align}
where the matrix $\vT_{(i)}$ extract an $m\times m$ frame  out of an $n\times n$ image. The matrix $T_{(i)}$ 
can be expressed in terms of a translation matrix $\vC_{\vx}$ which circularly translates  by $\vx$ in $D^{n\times n}$ and the
{\it restriction matrix} $\vR$, which is of size $m^2\times n^2$ so that $\vR(i,i)=1$ for all $i=1,\ldots,m^2$ and $0$ otherwise:
\begin{align}
\vT_{(i)}=&\vR\vC_{\vx_{(i)}}, \quad \vT=
\left (
\begin{array}{c}
 \vT_{(1)} \\
  \vT_{(2)} \\
  \vdots\\
\vT_{(K)} 
\end{array}
\right ) 
\end{align}


In the following, to simplify the notation, we will not distinguish between $w$ (resp. $\psi$) and its vector form $w^{\vee}$ (resp. $\psi^\vee$) and use the same notation $w$ (resp. $\psi$). We  simplify Eq (\ref{definition:Qi:rewrite}) as
\begin{align}
 \label{eq:standardprobe00} 
\vQ=\diag(\vS w)\vT.
\end{align}
Where we define $\vS\in \RR^{Km^2\times m^2}$ to be the matrix that replicates  the illumination probe $w$ into $k$ stack of frames, that is, $\vS$ is a $K\times 1$ block matrix with the $m^2\times m^2$ identity matrix as the block.
There is a  relationship between $\psi$,$w$ and $\vz$:
\begin{align}
 \label{eq:standardprobe0} 
 \vz=\vQ \psi=\diag(\vS w) \vT \psi=\diag(\vT \psi)  \vS w
\end{align}
 due to the fact that $\diag(\vS w) \vT \psi$ is the entry-wise product of $\vS w$ and $\vT\psi$.
%




In the notation used in this paper, projection operators are expressed as:
\begin{align*}
P_{\vQ}=\vQ(\vQ^*\vQ)^{-1}\vQ^*
\end{align*}
\begin{align*}
P_{\va}=\vF^{*}\frac{\vF\vz}{|\vF\vz|}\va
\end{align*}

\section{Alternating projections probe, frames, image \label{sec:frameprobeimage}}

The standard update\cite{Thibault2009, pie2}  proceeds in three steps updating the estimate of the illumination $w$,  estimate of the specimen $\psi$,
and the frames $\vz$ based on experimental geometry and data. First update the image $\psi$ by (\ref{eq:standardprobe0}):
\begin{align}
\psi\leftarrow &(\vQ^\ast \vQ)^{-1} \vQ^\ast \vz =\frac{\vT^\ast \diag (\vS \bar{w}) \vz}{\vT^\ast (\vS |w|^2) },
\end{align}
where the second equality holds due to (\ref{eq:standardprobe00}); that is, $\vQ^\ast \vQ=\diag(\vT^\ast (\vS |w|^2))$. Note that for the sake of simplifying notation, we represent $\diag(v)^{-1}B=:\frac{B}{v}$ when $\diag(v)$ is invertible diagonal matrix. Second, by (\ref{eq:standardprobe0}), update the probe $w$:
\begin{align}
\label{eq:standardprobe}
w\leftarrow&
\frac{\vS^\ast \diag( \vT \bar \psi ) \vz}{\vS^\ast \vT |\psi|^2 },
\end{align}
Third update the frames using the new probe embedded in $\vQ$, and the $P_{\va}$ operator:
\begin{align}
\label{eq:Pa}
\vz\leftarrow P_{\va}  \diag(\vS w) \vT \psi
\end{align} 
Note that Equation (\ref{eq:standardprobe}) is motivated by (\ref{eq:standardprobe0}). In fact, we have ${\diag(\vT \bar{\psi}) \vz=\diag(|\vT \psi|^2)}  \vS w$. 
And note that by a direct calculation, $\vS^*\diag(|\vT \psi|^2) \vS=\diag(\sum_i \vT_{(i)}|\psi|^2)$, which leads to  (\ref{eq:standardprobe}) when all entries of $\vS^\ast \vT |\psi|^2 $ are not zero. 
When there is zero entry in $\vS^\ast |\vT \psi|^2$, we may consider $\max\{\vS^\ast |\vT \psi|^2,\epsilon\}$, where $\epsilon>0$ is a small positive number determined by the user. Note that $\vS^*$ plays the role of averaging.

 \section{The relationship between probes (illumination ) and frames}
We now analyze here the symmetry between $\vQ$ and $\vz$, that is, we look at the  frame-wise relationship. 
We note that $T_{(i)}^\ast Q_{(i)}=T_{(i)}^\ast \diag (w) T_{(i)},\,\forall (i)$ is a diagonal matrix.
Hence we can write the following relationship between two frames:
\begin{align}
\nonumber
 T_{(i)}^\ast Q_{(i)} T_{(j)}^\ast z_{(j)}= T_{(i)}^\ast Q_{(i)} T_{(j)}^\ast Q_{(j)} \psi&= T_{(j)}^\ast Q_{(j)} T_{(i)}^\ast Q_{(i)} \psi
=   T_{(j)}^\ast  Q_{(j)} T_{(i)}^\ast  z_{(i)}
\nonumber
\end{align}

i.e. the i-th frame multiplied by the j-th illumination is equal to the j-th frame multiplied by the i-th illumination, after putting the results back to the right position.  
In other words, we have a pairwise relationship 
\begin{framed}
\begin{align}
\sum_{i,j} \|  T_{(j)}^\ast Q_{(j)} T_{(i)} ^\ast z_{(i)} -   T_{(i)} ^\ast Q_{(i)} T_{(j)}^\ast z_{(j)}\|=0.
\label{eq:pairwise}
\end{align}
\end{framed}

By swapping diagonal matrices $T_{(i)}^\ast Q_{(i)}$ 
and using $T_{(i)}T_{(i)}^\ast=I$, we can write out the pairwise discrepancy between all frames: 
\begin{align*}
\frac{1}{2}\sum_{i,j} \|  T_{(j)}^\ast Q_{(j)} T_{(i)} ^\ast z_{(i)} -   T_{(i)} ^\ast Q_{(i)} T_{(j)}^\ast z_{(j)}\|^2=&
\sum_{i,j} 
z^\ast_{(i)}  \left [   \mathbb{Q}_{(i)}^2  \delta_{i,j}
-   Q_{(i)} Q_{(j)}^\ast   \right ] z_{(j)}
=\vz^\ast [ \mathbb{Q}^2 - \mathbb{P}] \vz
\end{align*}
where ${\mathbb Q}^2_{(i)} \equiv  T_{(i)} \vQ^\ast \vQ  T_{(i)} ^\ast $, ${\mathbb Q}^2\equiv \diag({\mathbb Q}^2_{(i=1\dots k)})$,
$ \mathbb{P}\equiv \vQ \vQ^\ast  $. 
{
Indeed, note that $\sum_{ij}z_{(i)}^*\vT_{(i)}\vQ_{(j)}^*\vQ_{(j)} \vT_{(i)}^* z_{(i)}=\sum_{i}z_{(i)}^*\vT_{(i)}\big[\sum_{j}\vQ_{(j)}^*\vQ_{(j)}\big] \vT_{(i)}^* z_{(i)}$ and $\vQ^*\vQ=\sum_{j}\vQ_{(j)}^*\vQ_{(j)}$. Note that $\vQ^*\vQ$ is a diagonal matrix with non-negative diagonal entries describing how a pixel of the object of interest is illuminated by different windows. Thus we can define $\mathbb{Q}$ by taking square root of $\mathbb{Q}^2$ entry-wisely. We may assume that the diagonal is positive, otherwise some information of the object is missed in the experiment. Thus we may define $\mathbb{Q}^{-1}$. Also note that by a direct calculation, we obtain
\begin{align}\label{expansion:Q2}
\mathbb{Q}^2=\diag(\vT\vT^*\vS|w|^2).
\end{align}
 Intuitively, $\vT^*\vS|w|^2$ describes how pixels of the whole $\psi$ is illuminated, and $\vT_{(i)}\vT^*\vS|w|^2$ is how the pixels of the $i$-th patch of $\psi$ is illuminated, which is the same as $T_{(i)} \vQ^\ast \vQ  T_{(i)} ^\ast$.}

We can minimize this {functional} by first renormalizing, making a change of variable $\hat {\vz}= \mathbb{Q}\vz$, then applying the power iteration:
\begin{align}
\label{eq:Q2P}
\vz^\ast [ \mathbb{Q}^2 - \mathbb{P}] \vz
=\hat{\vz}^\ast [ I -\mathbb{Q}^{-1} \mathbb{P} \mathbb{Q}^{-1}] \hat \vz
=\hat{\vz}^\ast [ I -P_Q] \hat \vz\\
\vz\leftarrow \mathbb{Q}^{-1} P_Q  {\hat \vz}=\mathbb{Q}^{-1} P_Q \mathbb{Q}  {\vz}=P_Q  {\vz}
\label{eq:projection_update}
\end{align}
which holds by the way $\mathbb{Q}$ is defined. 
The reason for the above formulation is to establish a relationship between $Q$ and $z$.
As we have seen in Sec. \ref{sec:frameprobeimage}, the  relationship between $w$ and $\psi$ in (Eq. (\ref{eq:standardprobe}))
 is used to recover both $w$ and $\psi$. 
 
Interestingly, the relationship (Eq. \ref{eq:pairwise}) is symmetric w.r.t. $Q$ (i.e. $w$) and $z$. In other words, we can update the probe based on the frames $z$
by solving the symmetric counterpart: 
\begin{align*}
\frac{1}{2}\sum_{i,j} \|  T_{(j)}^\ast Q_{(j)} T_{(i)} ^\ast z_{(i)} -   T_{(i)} ^\ast Q_{(i)} T_{(j)}^\ast z_{(j)}\|^2=&
\vz^\ast [ \mathbb{Q}^2 - \mathbb{P}] \vz=\vz^\ast [\diag( \vT \vT^\ast \vS |w|^2) - \diag(\vS w) \vT \vT^\ast \diag(\vS\bar w)] \vz\\
=&w^\ast \vS^\ast [\diag(\vT\vT^\ast |\vz|^2) -  \diag(\vz) \vT\vT^\ast\diag(\bar \vz)] \vS w,
\end{align*}
{where the second equality holds due to (\ref{expansion:Q2}) and (\ref{eq:standardprobe00}), and the last equality holds due to the equality $\vz^*\diag( \vT \vT^\ast \vS |w|^2)\vz=w^\ast \vS^\ast \diag(\vT\vT^\ast |\vz|^2)\vS w$ by a direct calculation. Note that 
\begin{enumerate}
\item we may view $\vz^*\diag( \vT \vT^\ast \vS |w|^2)\vz$ as a weighted inner product in $\CC^{Km^2}$;
\item when $z\in\CC^{Km^2}$ is the vectorized version of all illuminated images, geometrically $\vS^*\diag(z)\vS$ means averaging over all illuminated images; that is $\vS^*\diag(z)\vS=\diag(\sum_{i=1}^K{z_{(i)}})$. By a directly calculation, we have $\vS^*\diag(z)\vS= \diag(\vS^*z) $.
\end{enumerate}
}
{With the above preparation, we} wish to solve:
\begin{align}\label{probe_newmin}
\arg\min_{w} w^\ast \vS^\ast  [  \diag ( \vT \vT^\ast  | \vz|^2) - \diag( \bar \vz^\ast )\vT \vT^\ast \diag( \bar \vz ) ]\vS w,  
\end{align} 
which can be expressed as (in a similar form as Eq.  (\ref{eq:Q2P}))
\begin{align}\label{eq:DA}
&\arg\min_w w^\ast [D-A] w,\\ \mbox{where }
&\left\{
\begin{array}{rl}
D=&\vS^\ast  \diag ( \vT \vT^\ast  | \vz|^2) \vS=  \diag(\vS^\ast  ( \vT \vT^\ast  | \vz|^2)) \\
A=&\vS^\ast  [   \diag( \bar \vz^\ast )\vT \vT^\ast \diag( \bar \vz ) ]\vS  .
\end{array}
\right. \nonumber
\end{align}
Note that $A$ is Hermitian but not a diagonal matrix since $\diag( \bar \vz^\ast )\vT \vT^\ast \diag( \bar \vz )$ is not diagonal.
Thus, (\ref{eq:DA}) yields - by power method, starting from an initial estimate $w$:
\begin{align}\label{newProbeUpdate}
w \leftarrow D^{-1} Aw 
=\frac{\vS^\ast \diag( \bar \vz^\ast )\vT \vT^\ast \diag( \bar \vz ) \vS w}{\vS^\ast   \vT \vT^\ast  | \vz|^2}
=\frac{\vS^\ast \diag( \vz )\vT \vT^\ast \diag( \vS w  ) \bar \vz }{\vS^\ast   \vT \vT^\ast  | \vz|^2}
=\frac{\vS^\ast \diag(  \vz )\vT  \vQ^\ast  \bar \vz }{\vS^\ast   \vT \vT^\ast  | \vz|^2},
\end{align}
where the last equality holds due to (\ref{eq:standardprobe00}).
 Note the projection  update (\ref{eq:projection_update}) can be expressed in a similar way 
by swapping the order of $\vS$ with $\vT \vT^\ast$:
$$
\vz\leftarrow P_{\vQ} \vz=
\frac{\vQ\vQ^\ast \vz}{\QQ^2} =
 \frac{\diag (\vS w) \vT \vT^\ast \diag(  \vS w)^\ast \vz
 }{\vT \vT^\ast \vS |w|^2} 
$$


How does this update relate to the standard update {in (\ref{eq:standardprobe})} {{blue}If} we insert $\QQ^2$ inside the averaging matrices $\vS^\ast () \vS$ and when $z=\vQ\psi$, we obtain
$$
w=\frac{\vS^\ast {\mathbb{Q}^{-2}} \diag(  \vz )\vT \vT^\ast \diag( \bar \vz ) \vS w}{\vS^\ast  {\mathbb{Q}^{-2}}   \vT \vT^\ast  | \vz|^2}
=\frac{\vS^\ast \diag(  \vz )\vT \left[\frac{1}{\vQ^\ast \vQ} \vT^\ast \diag( \bar \vz ) \vS w\right]}{\vS^\ast   \vT \left[\frac{1}{\vQ^\ast \vQ}\vT^\ast  | \vz|^2\right]}
= \frac{\vS^\ast \diag(\vT\bar \psi)  \vz}{\vS^\ast   \vT | \psi'|^2};
$$
{where the third equality holds since $\vT^\ast \diag( \bar \vz ) \vS w=\vT^\ast \diag( \vS w)\bar \vz =\overline{\vT^\ast \diag( \vS \bar w) \vz }= \overline{\vQ^* z}$ and $\psi=(\vQ^*\vQ)^{-1}\vQ^*z$, and we define $|\psi' |^2 \equiv  \frac{1}{\vQ^\ast \vQ} \vT^\ast  | \vz|^2 $. }
Note that $\frac{1}{\vQ^\ast \vQ} \vT^\ast  | \vz|^2 $ is different from $|\psi|^2=|\frac{1}{\vQ^\ast \vQ}\vQ^\ast \vz|^2$ {in (\ref{eq:standardprobe})}. However  $\vS^\ast \vT$ in (\ref{newProbeUpdate}) smears out the normalization factor, and the two results are similar.
  

\section{Rank-1 speedup for weakly scattering and piecewise smooth objects}

A simple way to speed up is to simply remove  one constant term (DC) on the fly, that is we estimate the  \textit{transparency} factor $\nu\in \CC^1$
$$
 \nu= \frac{\| \vQ^\ast \vQ \psi \|}{\| \vQ^\ast \vQ \|}=
 \frac{(\vS w)^\ast \vz}{\|\vS w\|^2}
$$
where $\vz^{\ell}=\vz^{}$
and subtract the average transmitted intensity $\vz\leftarrow\vz-\nu \vS w$ in Eqs. (\ref{eq:DA}, \ref{newProbeUpdate}) 

we get the following update:
\begin{framed}


\begin{align}
\label{eq:probeupdate2}
,\quad w
\leftarrow\frac{\vS^\ast\left (  \diag(  \vz -\nu \vS w) \vT \bar{\vQ}^\ast (\bar{z}-\bar \nu \vS \bar w)  \right ) }{\vS^\ast  ( \vT \vT^\ast  | \vz-\nu \vS w|^2  )}
=\frac{\vS^\ast\left (  \diag(  \vz -  \nu w) ( \vT \bar \vQ^\ast \bar z- \bar \nu  \QQ^2 \right ) }{\vS^\ast  ( \vT \vT^\ast  | \vz|^2 -2 \Re (\bar \nu \vT \vQ^\ast \vz)+\QQ^2 |\nu|^2 )},
\end{align}
where the last equality holds by (\ref{expansion:Q2}).
%


\end{framed}

{
This is useful when an object has a  strong DC term (weak contrast).  Moreover, we can remove the  frame-wise DC term for piecewise objects, which is useful since any constant region within the object does not provide any information about the probe.  The second formulation enables us to compute $\nu$ frame-wise. 


\begin{framed}

We simply average the frames that overlap together, that is we apply the  operator that sums all the frames that overlap with a given frame. 
Consider the $K\times K$ matrix $\vX$, with entries $\vX{(i,j)}=1$ if $(i)$ overlaps with $(j)$, or 0 otherwise.  Compute the $\upsilon\in\CC^K$ vector:
\begin{align*}
\upsilon_{i}= \frac{\sum_{j} \vX_{i,j}  ( w^\ast \vz_{(j)} )}{\sum_{j} \vX_{i,j}  (\|w \|^2 )}
= \frac{\sum_{j} \vX_{i,j}  ( w^\ast \vz_{(j)} )}{ \|w \|^2  \sum_{j} \vX_{i,j} }
\end{align*}
Then do the following update, which we write with some abuse of notation:
\begin{align*}
w&\leftarrow\frac{\vS^\ast\left (  \diag(  \vz -  \upsilon \vS w) ( \vT \bar \vQ^\ast \bar z- \bar \upsilon  \QQ^2 \right ) }{\vS^\ast  ( \vT \vT^\ast  | \vz|^2 -2 \Re (\bar \upsilon \vT \vQ^\ast \vz)+\QQ^2 |\upsilon|^2 )}
\end{align*}
the abuse of notation is that  $(\vS w) \upsilon$ , $\QQ \bar \upsilon$ and others, are intended as   $(\vS w) \vB \upsilon$ where $\vB$ is the matrix that replicates $\upsilon_{i}$ onto the frame of dimension $n^2$.
\end{framed}
}

Numerical tests are shown in Fig. \ref{fig:plot} with $n=223$, $m=128$,  $K=20\times20$, 5 pixels steps.
Significant speedup is observer using the update in Eq. (\ref{eq:probeupdate2}), vs Eq. (\ref{eq:standardprobe}), 
the speedup is accentuated when there is a strong average constant transmission factor.

\begin{figure}[hbtp]
  \begin{center}
  \subfigure[True illumination]{
     \includegraphics[width=0.3\textwidth]{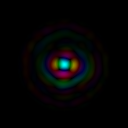}
}
  \subfigure[Initial estimate of the illumination]{
     \includegraphics[width=0.3\textwidth]{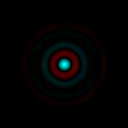}
}
\\
  \subfigure[Final estimate of the illumination]{
     \includegraphics[width=0.3\textwidth]{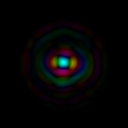}
}
  \subfigure[NRMSE $\min_{c^{\ell}} \|c^{\ell}w^{(\ell)}-w_{0}\|/\|w_{0}\|$]{
               \includegraphics[width=0.3\textwidth,clip]{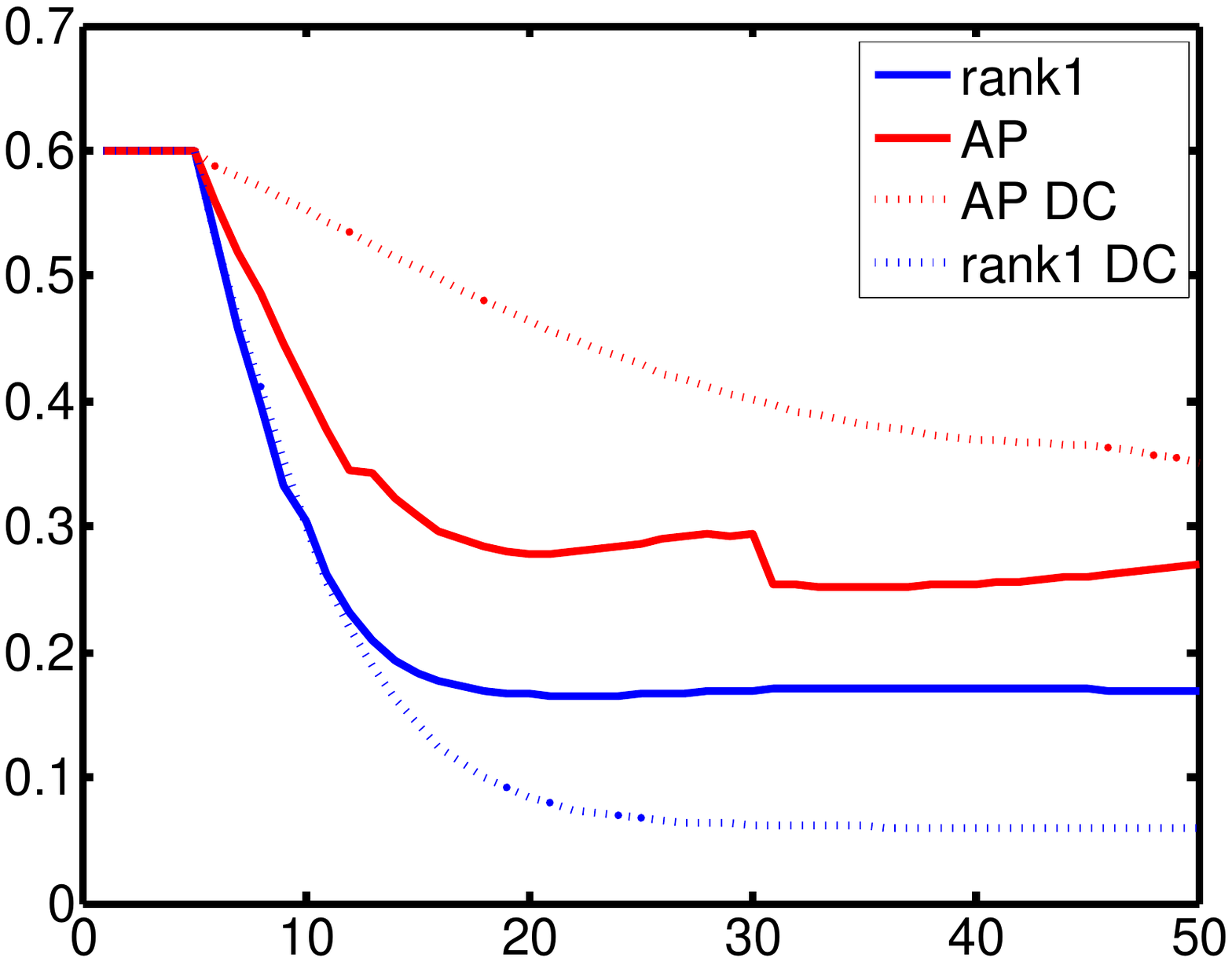}  
}\\
  \subfigure[True specimen]{
                              \includegraphics[width=0.3\textwidth]{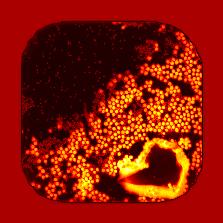}        
}
  \subfigure[Recovered specimen]{
                              \includegraphics[width=0.3\textwidth]{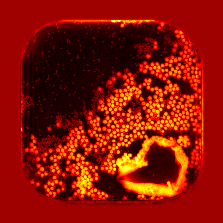}        
}
 \end{center}
 \caption{
\label {fig:plot}
Initial probe, true illumination, final illumination.
Convergence with strong DC term (99\%) or with DC term removed.
test object, reconstruction. Residual disparities are ascribed to global translations 
between the ``true'' illumination and image used for simulations. At every illumination update, the center of mass of the illuminations is  centered within one neighboring pixel. Further analysis of fractional shifts and aberrations can be performed using the prescription given in \cite{marchesiniXRM}.
}
\end{figure}

\section{Acknowledgements}
This work is partially supported by the Center for Applied Mathematics for Energy Research Applications (CAMERA), which is a partnership between Basic Energy Sciences (BES) and  Advanced Scientific Computing Research (ASRC)  at the U.S Department of Energy under contract DE-AC02-05CH11231 (SM).
\hrulefill

\noindent

\end{widetext}
\bibliographystyle{plain}
\bibliography{ptycho_probe}
\end{document}